\documentclass{amsart}
\usepackage{epsfig}
\usepackage{amssymb}

\newtheorem{theorem}{Theorem}[section]
\newtheorem{prop}[theorem]{Proposition}
\newtheorem{lemma}[theorem]{Lemma}
\newtheorem{corollary}[theorem]{Corollary}
\newcommand\beq{\begin{equation}}
\newcommand\eeq{\end{equation}}
\newcommand\bce{\begin{center}}
\newcommand\ece{\end{center}}
\newcommand\bea{\begin{eqnarray}}
\newcommand\eea{\end{eqnarray}}
\newcommand\ben{\begin{enumerate}}
\newcommand\een{\end{enumerate}}
\newcommand\brr{\begin{array}}
\newcommand\err{\end{array}}
\newcommand\bt{\begin{tabular}}
\newcommand\et{\end{tabular}}
\newcommand\nin{\noindent}
\newcommand\nn{\nonumber}

\newcommand\ms{\medskip}

\newcommand\ol{\overline}
\renewcommand\S{{\mathcal S}}
\def\mn{\mbox{-}}

\def\T{{\mathcal T}}
\def\H{\tilde{H}}

\title{Generating trees for permutations avoiding generalized patterns}
\author{Sergi Elizalde}
\address{Department of Mathematics, Dartmouth College, Hanover, NH
03755 / Centre de Recerca Matem\`atica, E-08193 Bellaterra, Spain}
\email{sergi.elizalde@dartmouth.edu}

\begin{document}
\maketitle \vspace{-8mm}
\begin{abstract}
We construct generating trees with one, two, and three labels for
some classes of permutations avoiding generalized patterns of length
$3$ and $4$. These trees are built by adding at each level an entry
to the right end of the permutation, which allows us to incorporate
the adjacency condition about some entries in an occurrence of a
generalized pattern. We use these trees to find functional equations
for the generating functions enumerating these classes of
permutations with respect to different parameters. In several cases
we solve them using the kernel method and some ideas of
Bousquet-M\'elou \cite{B-M}. We obtain refinements of known
enumerative results and find new ones.
\end{abstract}

\section{Introduction}

\subsection{Generalized pattern avoidance}
We denote by $\S_n$ the symmetric group on $\{1,2,\ldots,n\}$. Let
$n$ and $k$ be two positive integers with $k\le n$, and let
$\pi=\pi_1\pi_2\cdots\pi_n\in\S_n$ be a permutation. A generalized
pattern $\sigma$ is obtained from a permutation
$\sigma_1\sigma_2\cdots\sigma_k\in\S_k$ by choosing, for each
$j=1,\ldots,k-1$, either to insert a dash $\mn$ between $\sigma_j$
and $\sigma_{j+1}$ or not. More formally,
$\sigma=\sigma_1\varepsilon_1\sigma_2\varepsilon_2\cdots\varepsilon_{k-1}\sigma_k$,
where each $\varepsilon_j$ is either the symbol~$\mn$ or the empty
string. With this notation, we say that $\pi$ \emph{contains} (the
generalized pattern) $\sigma$ if there exist indices
$i_1<i_2<\cdots<i_k$ such that (i) for each $j=1,\ldots,k-1$, if
$\varepsilon_j$ is empty then $i_{j+1}=i_j+1$, and (ii) for every
$a,b\in\{1,2,\ldots,k\}$, $\pi_{i_a}<\pi_{i_b}$ if and only if
$\sigma_a<\sigma_b$. In this case,
$\pi_{i_1}\pi_{i_2}\cdots\pi_{i_k}$ is called an \emph{occurrence}
of $\sigma$ in $\pi$.

If $\pi$ does not contain $\sigma$, we say that $\pi$ \emph{avoids}
$\sigma$, or that it is \emph{$\sigma$-avoiding}. For example, the
permutation $\pi=3542716$ contains the pattern $12\mn4\mn3$ because
it has the subsequence $3576$. On the other hand, $\pi$ avoids the
pattern $12\mn43$. We denote by $\S_n(\sigma)$ the set of
permutations in $\S_n$ that avoid $\sigma$. More generally, if
$\Sigma=\{\sigma_1,\sigma_2,\ldots\}$ is a collection of generalized
patterns, we say that a permutation $\pi$ is $\Sigma$-avoiding if
$\pi$ is $\sigma$-avoiding for all $\sigma \in \Sigma$. We denote by
$\S_n(\Sigma)$ the set of $\Sigma$-avoiding permutations in~$\S_n$.

We use the word {\em length} to refer to the number of letters in a
permutation, so that $\S_n$ is the set of permutations of length
$n$. A {\em class} will consist of a set (e.g., all permutations
avoiding a given pattern) together with a function (e.g., the
length). Given a permutation $\pi\in\S_n$, we will write
$r(\pi)=\pi_n$ to denote the rightmost entry of $\pi$. In all our
generating functions, the variable $t$ will mark the length of the
permutation.

\subsection{Generating trees}

Generating trees are a useful tool for enumerating classes of
pattern-avoiding permutations (see, for example, \cite{West,Wes96}).
The nodes at each level of the generating tree are indexed by
permutations of a given length. It is common in the literature to
define the children of a permutation $\pi$ of length $n$ to be those
permutations that are obtained by inserting the entry $n+1$ in
$\pi=\pi_1\pi_2\cdots\pi_n$ in such a way that the new permutation
is still in the class. In this paper we consider a variation of this
definition. Here, the children of a permutation $\pi$ of length $n$
are obtained by appending an entry to the right of $\pi$, and adding
one to all the entries in $\pi$ that were greater than or equal to
the new entry. For example, if the entry $3$ is appended to the
right of $\pi=24135$, the child that we obtain is $251463$. Adding
the new entry to the right of the permutation makes these trees
well-suited to enumerate permutations avoiding generalized patterns,
as we will see throughout the paper. We will refer to these trees as
{\em rightward generating trees}. This kind of generating trees has
been used in~\cite{BFP} to enumerate permutations avoiding sets of
three generalized patterns of length three with one dash, such as
$\{1\mn23,2\mn13,1\mn32\}$.

For some classes of permutations, a label $(\ell)$ can be associated
to each node of the tree in such a way that the number of children
of a permutation and their labels depend only on the label of the
parent. For example, in the tree for $1\mn2\mn3$-avoiding
permutations, we can label each node $\pi$ with
$m=\min\{\pi_i:\exists j<i \mbox{ with } \pi_{j}<\pi_i\}$ (or
$m=n+1$ if $\pi=n\cdots21$). Then, the children of a permutation
with label $(m)$ have labels $(m+1),(2),(3),\ldots,(m)$,
corresponding to the appended entry being $1,2,3,\ldots,m$,
respectively. This succession rule, together with the fact that the
root ($\pi=1\in\S_1$) has label $(2)$, completely determines the
tree. From this rule one can derive a functional equation for the
generating function that enumerates the permutations by their length
and the label of the corresponding node in the tree for this class
of $1\mn2\mn3$-avoiding permutations. For generating trees with one
label, these equations are well understood and their solutions are
algebraic series. This is the case of the generating trees obtained
in~\cite{BFP}, for example.

In other cases, however, one label is not enough to describe the
generating tree in terms of a succession rule.
 Generating trees with
two labels were used in \cite{B-M} to enumerate restricted
permutations. In fact, the inspiration for the present paper and
many of the ideas used come from Bousquet-M\'elou's work. One
difference is that here trees are constructed by adding at each
level an entry to the right end of the permutation, which allows us
to keep track of elements occurring in adjacent positions. In
Section~\ref{sec:threelab} we consider some classes of permutations
whose rightward generating tree has three labels for each node.

\subsection{Organization of the paper}

In this paper we enumerate several families of permutations that
avoid generalized patterns. What ties together the results in the
different sections is the technique that we use to obtain them. The
strategy consists of building a rightward generating tree for the
family of permutations, translating the succession rule into a set
of functional equations,  and applying the kernel method to them. We
have tried this strategy for a number of classes of permutations,
and we have found it to work in several cases, which we include
here. This is why the sets of generalized patterns that we discuss
may seem somewhat arbitrary. For other patterns one can construct
similar generating trees with two or three labels, but we have not
been able to solve the corresponding functional equations for the
generating function,  so we have not included these examples here.
In any event, this paper is not meant to be an exhaustive study of
the sets of patterns for which this technique would work.

In general, we have looked for sets of patterns for which the
rightward generating tree of the class of permutations avoiding them
has a simple succession rule, once appropriate labels are chosen. In
some cases, we have chosen patterns based on the elegance of their
enumerating sequence, like in Section~\ref{sec:oddbar}. In others, we have
chosen patterns whose corresponding generating function has zero
radius of convergence, as is the case in Sections
\ref{sec:1_23_3_12}, \ref{sec:threelabsub1}, and
\ref{sec:threelabsub2}. These seem to be the first instances of
generating functions with zero radius of convergence that arise from
generating trees and the kernel method.

We have classified the sets of studied patterns depending on how
many labels are needed to describe the generating tree. In
Section~\ref{sec:onelab} we consider some families of permutations
where the tree can be described with one label, which is the value
of the rightmost entry in the permutation. The results in this
section are new, and all involve permutations that avoid the pattern
$2\mn1\mn3$. This makes the succession rules easier because this
restriction prevents a permutation with rightmost value $r$ to have
a child with rightmost value greater than $r+1$.

In Section~\ref{sec:twolab} we study classes of permutations where
each node of the generating tree bears a pair of labels. For most of
them we get rational or algebraic generating functions, and their
enumeration has been done in the literature using different
techniques.
 Section~\ref{sec:threelab} contains some of the main results of the
paper. We find ordinary generating functions for $\{1\mn23,3\mn12,34\mn21\}$-avoiding and
$\{1\mn23,34\mn21\}$-avoiding permutations. Both families are
described by generating trees with three labels.

Additional motivation for the study of these families of
permutations comes from trying to understand the possible asymptotic
behaviors of the number of permutations avoiding generalized
patterns (see~\cite{Eli06}). An asymptotic analysis of the
coefficients of the generating functions for
$\{1\mn23,3\mn12\}$-avoiding, $\{1\mn23,34\mn21\}$-avoiding, and
$\{1\mn23,3\mn12,34\mn21\}$-avoiding permutations that we have found
may reveal that their asymptotic growth is strictly smaller than
that of Bell numbers but strictly greater than exponential. This
would be the first known instance of a family of pattern-avoiding
permutations that exhibits such a behavior.

\section{Generating trees with one label}\label{sec:onelab}

In this section we enumerate classes of pattern-avoiding
permutations whose rightward generating trees can be described by a
succession rule involving only one label for each node. For some of
these classes, rightward generating trees are not the only way to
obtain the results, but they are a tool that works in all these
cases.

The classes in this section avoid the pattern $2\mn1\mn3$. Note that
avoiding this pattern is equivalent to avoiding the generalized
pattern $2\mn13$. Indeed, if $\pi$ contains an occurrence of
$2\mn1\mn3$, say $\pi_i\pi_j\pi_k$  with $\pi_j<\pi_i<\pi_k$, then
there must be some index $\ell$ with $j\le \ell<k$ such that
$\pi_\ell<\pi_i$ and $\pi_{\ell+1}>\pi_i$, so
$\pi_i\pi_\ell\pi_{\ell+1}$ is an occurrence of $2\mn13$. For any
class of permutations that avoid this pattern, the corresponding
rightward generating tree has the property that the appended entry
at each level can never be more than one unit larger than the entry
appended at the previous level.

\subsection{$\{2\mn1\mn3,\ol{2}\mn31\}$-avoiding permutations}
\label{sec:213bar231}

A permutation $\pi$ is said to avoid the \emph{barred} pattern
$\ol{2}\mn31$ if every descent in $\pi$ (an occurrence of the
generalized pattern $21$) is part of an occurrence of $2\mn31$;
equivalently, for any index $i$ such that $\pi_i>\pi_{i+1}$ there is
an index $j<i$ such that $\pi_i>\pi_j>\pi_{i+1}$. The bar indicates
that the $2$ is forced whenever a $31$ occurs. For example, the
permutation $4627513$ avoids $\ol{2}\mn31$, but $2475613$ does not.

We use $M_n$ to denote the $n$-th Motzkin number. Recall that
$\sum_{n\ge0} M_{n} t^n=\frac{1-t-\sqrt{1-2t-3t^2}}{2t^2}$. The next
result seems to be a new interpretation of the Motzkin numbers.

\begin{prop}\label{th:213bar231}
The number of $\{2\mn1\mn3,\ol{2}\mn31\}$-avoiding permutations of
size $n$ is $M_{n-1}$.
\end{prop}

\begin{proof}
Consider the rightward generating tree for
$\{2\mn1\mn3,\ol{2}\mn31\}$-avoiding permutations. Labeling each
permutation with its rightmost entry $r=r(\pi)$, this tree is
described by the succession rule \bce\bt{l}
$(1)$ \\
$(r)\longrightarrow (1)\ (2)\ \cdots\ (r-1)\ (r+1)$. \et\ece
 Indeed, the new entry appended to the right of $\pi$ cannot be greater
than $\pi_n+1$ in order for the new permutation to be
$2\mn1\mn3$-avoiding, and it cannot be $\pi_n$ because then it would
create an occurrence of $21$ that is not part of an occurrence of
$2\mn31$.

Defining $D(t,u)=\sum_{n\ge1}\
\sum_{\pi\in\S_n(2\mn1\mn3,\ol{2}\mn31)} u^{r(\pi)} t^n$,
 the succession rule above gives the following equation for the
generating function:
\beq\label{eq:kernel_D}\left(1-\frac{t}{u-1}-tu\right)D(t,u)=tu-\frac{tu}{u-1}D(t,1).\eeq
The next step is to apply the kernel method. This technique, which
has been part of mathematical folklore for decades, has recently
been systematized in~\cite{BBDFGG,BF,BP}. Of the two values of $u$
as a function of $t$ that cancel the term multiplying $D(t,u)$ on
the left hand side, $u_0=u_0(t)=\frac{1+t-\sqrt{1-2t-3t^2}}{2t}$ is
a well-defined formal power series in $t$. Substituting $u=u_0$ in
(\ref{eq:kernel_D}) gives
$$D(t,1)=u_0-1=\frac{1-t-\sqrt{1-2t-3t^2}}{2t},$$ which is the
generating function for the Motzkin numbers with the indices shifted
by one.
\end{proof}

There is also a bijective proof of Proposition~\ref{th:213bar231}.
Given a permutation
$\pi=\pi_1\pi_2\cdots\pi_n\in\S_n(2\mn1\mn3,\ol{2}\mn31)$, we can
construct a Dyck path of size $n$ (i.e., a sequence of $n$ $U$s and
$n$ $D$s so that no prefix contains more $D$s than $U$s) as follows.
A right-to-left maximum of $\pi$ is an entry $\pi_i$ such that
$\pi_i>\pi_j$ for all $j>i$. Let
$\pi_{i_1},\pi_{i_2},\ldots,\pi_{i_m}$ be the right-to-left maxima
of $\pi$, with $i_1<i_2<\cdots<i_m=n$. Consider the Dyck path
$$\varphi(\pi)=U^{i_1}D^{\pi_{i_1}-\pi_{i_{2}}}U^{i_2-i_1}D^{\pi_{i_2}-\pi_{i_{3}}}U^{i_3-i_2}\cdots
D^{\pi_{i_{m-1}}-\pi_{i_{m}}}U^{i_m-i_{m-1}}D^{\pi_{i_m}},$$ where
exponentiation indicates repetition of a step. This map is a
bijection between $2\mn1\mn3$-avoiding permutations and Dyck paths
(see~\cite{Kra}), and it is not hard to see that the condition if
$\pi$ being $\ol{2}\mn31$-avoiding is equivalent to the requirement
that the path contains no three consecutive steps $UDU$. So, we have
a bijection between $\S_n(2\mn1\mn3,\ol{2}\mn31)$ and $UDU$-free
Dyck paths of size $n$.

To finish the proof, we next describe a bijection due to
Callan~\cite{Ca04} between $UDU$-free Dyck paths of size $n$ and
Motzkin paths of length $n-1$ (i.e., sequences of $n-1$ steps $U$,
$D$, and $H$ with the same number of $U$s and $D$s and so that no
prefix contains more $D$s than $U$s). We say that a $U$ and a $D$ in
a Dyck path are matched if the $D$ is to the right of the $U$ and
the letters between them form a Dyck path. Note that each step is
matched with exactly another one. Given a $UDU$-free Dyck path,
first append a $D$ to it. Now, for each $D$ that is immediately
preceded and followed by $D$ steps, delete it and replace its
matching $U$ with an $H$. Next, replace each occurrence of $UDD$
with a $D$. Finally, delete the $D$ that was appended to the path.
This produces a Motzkin path of length $n-1$. The composition of
these two bijections completes the bijective proof of
Proposition~\ref{th:213bar231}.

\subsection{$\{2\mn1\mn3,\ol{2}^o\mn31\}$-avoiding
permutations}\label{sec:oddbar}

Extending the notion of barred patterns, we say that a permutation
$\pi$ avoids the pattern $\ol{2}^o\mn31$ if every descent in $\pi$
is the `31' part of an odd number of occurrences of $2\mn31$;
equivalently, for any index $i$ such that $\pi_i>\pi_{i+1}$, the
number of indices $j<i$ such that $\pi_i>\pi_j>\pi_{i+1}$ is odd.

\begin{prop}\label{th:odd}
The number of $\{2\mn1\mn3,\ol{2}^o\mn31\}$-avoiding permutations of
size $n$ is
\beq\label{eq:cat3} |\S_n(2\mn1\mn3,\ol{2}^o\mn31)|=\begin{cases} \frac{1}{2k+1}\binom{3k}{k} & \mbox{if } n=2k, \\
\frac{1}{2k+1}\binom{3k+1}{k+1} & \mbox{if } n=2k+1. \end{cases}\eeq
\end{prop}

\begin{proof}
The rightward generating tree for
$\{2\mn1\mn3,\ol{2}^o\mn31\}$-avoiding permutations is given by the
succession rule \bce\bt{l}
$(1)$ \\
$(r)\longrightarrow  \cdots(r-3)\ (r-1)\ (r+1) $, \et\ece that is,
the labels of the children of a node labeled $r$ are the numbers
$1\le j\le r+1$ such that $r-j$ is odd.  Let $J(t,u)=\sum_{n\ge1}\
\sum_{\pi\in\S_n(2\mn1\mn3,\ol{2}^o\mn31)} u^{r(\pi)}
t^n=\sum_{r\ge1} J_r(t) u^r$, and let $J^e(t,u)=\sum_{r\
\mathrm{even}} J_r(t) u^r$. The succession rule translates into the
following functional equation:
\beq\label{eq:kernel_J}\left(1-\frac{tu^3}{u^2-1}\right)J(t,u)=tu-\frac{tu^2}{u^2-1}J(t,1)+\frac{tu(u-1)}{u^2-1}J^e(t,1).\eeq
The kernel $1-\frac{tu^3}{u^2-1}$ as a function in the variable $u$
has three zeroes, two of which are complex conjugates. Denote them
by $u_1=a(t)+b(t)i$ and $u_2=\bar{u_1}=a(t)-b(t)i$. Adding the
equations $0=u_i^2-1-u_i J(t,1)+(u_i-1)J^e(t,1)$ for $i=1,2$, we get
$$a(t) J(t,1)=a(t)^2-b(t)^2-1+(a(t)-1)J^e(t,1),$$ and subtracting them
gives $$J(t,1)=2a(t)+J^e(t,1).$$ Solving this system of equations
for $J$, we get that $J(t,1)=2a(t)-a(t)^2-b(t)^2-1$. Plugging in the
values of $a(t)$ and $b(t)$ yields the expression
$$J(t,1)=\frac{(2-3t)f(t)^2+(9t-2-g(t))f(t)+(2-6t)g(t)+54t^2-18t-4}{3tf(t)^2},$$
where $g(t)=\sqrt{3(27t^2-4)}$ and $f(t)=[12tg(t)-108t^2+8]^{1/3}$.
It is easy to check that $J=J(t,1)$ is a root of the polynomial
$tJ^3+(3t-2)J^2+(3t-1)J+t=0$. Using the Lagrange inversion formula,
one sees that its coefficients are given by (\ref{eq:cat3}), which
is sequence A047749 from the On-Line Encyclopedia of Integer
Sequences~\cite{Slo}. Observe that we can also obtain an expression
for $J(t,u)$ using (\ref{eq:kernel_J}) and the fact that
$J^e(t,1)=-a(t)^2-b(t)^2-1$.
\end{proof}

It is also possible to give a direct bijective proof of
Proposition~\ref{th:odd} that does not use rightward generating
trees. A well-known combinatorial interpretation of the numbers
(\ref{eq:cat3}) is that they enumerate lattice paths from $(0,0)$ to
$(n,\lfloor n/2\rfloor)$ with steps $E=(1,0)$ and $N=(0,1)$ that
never go above the line $y=x/2$. We next describe a bijection from
$\S_n(2\mn1\mn3,\ol{2}^o\mn31)$ to these paths. Let
$\pi=\pi_1\pi_2\cdots\pi_n\in\S_n(2\mn1\mn3,\ol{2}^o\mn31)$. Let
$\pi_{i_1},\pi_{i_2},\ldots,\pi_{i_m}$ be the right-to-left maxima
of $\pi$, with $i_1<i_2<\cdots<i_m=n$. We claim that the condition
that $\pi$ is $\{2\mn1\mn3,\ol{2}^o\mn31\}$-avoiding guarantees that
all the differences $\pi_{i_j}-\pi_{i_{j+1}}$ are even. To see this,
fix $j$ and let $\mathcal{O}$ be the set of entries $a$ such that
$a\pi_{i_j}\pi_{i_j+1}$ is an occurrence of $2\mn31$. Since $\pi$
avoids $\ol{2}^o\mn31$, the cardinality of $\mathcal{O}$ is odd.
Now, every $a\in\mathcal{O}$ must satisfy $a>\pi_{i_{j+1}}$. This is
obvious if $i_j+1=i_{j+1}$, and otherwise it follows from the fact
that if $a<\pi_{i_{j+1}}$, then $a\pi_{i_j+1}\pi_{i_{j+1}}$ would be
an occurrence of $2\mn1\mn3$. On the other hand, any entry $a$ with
$\pi_{i_{j+1}}<a<\pi_{i_j}$ must appear to the left of $\pi_{i_j}$,
since $\pi_{i_j}$ and $\pi_{i_{j+1}}$ are consecutive right-to-left
maxima, and so $a\in\mathcal{O}$. Thus, $\mathcal{O}$ is precisely
the set of integers strictly between $\pi_{i_j}$ and
$\pi_{i_{j+1}}$, which implies that $\pi_{i_j}-\pi_{i_{j+1}}$ is
even. Now, for $j=1,\ldots,m-1$, let
$a_j=(\pi_{i_j}-\pi_{i_{j+1}})/2$. Let
$a_m=\lfloor\pi_{i_m}/2\rfloor$. We map $\pi$ to the following path
from $(0,0)$ to $(n,\lfloor n/2\rfloor)$:
$$E^{i_1}N^{a_1}E^{i_2-i_1}N^{a_2}E^{i_3-i_2}N^{a_3}\cdots
E^{i_m-i_{m-1}}N^{a_m}.$$ It can be checked that this map is a
bijection.
 For example, if $\pi=4675123$, we have
$\pi_{i_1}\pi_{i_2}\pi_{i_3}=\pi_{3}\pi_{4}\pi_{7}=753$, so the
corresponding path from $(0,0)$ to $(7,3)$ is $EEENENEEEN$.

Aside from lattice paths, the sequence
$d_n:=|\S_n(2\mn1\mn3,\ol{2}^o\mn31)|$ from (\ref{eq:cat3}) is also
known to enumerate symmetric ternary trees with $3n$ edges and
symmetric diagonally convex directed polyominoes of area $n$. These
numbers have also appeared before in connection to pattern-avoiding
permutations. It is shown in~\cite{BEM} that the number of
$2143$-avoiding Dumont permutations of the second kind of length
$2n$ is $d_n d_{n+1}$ (see~\cite{BEM} for definitions). The sequence
$d_n$ enumerates what the authors call {\it lower boards}, which are
$2\mn1\mn3$-avoiding permutations of length $n$ whose diagram fits
in a certain shape. A bijection between such permutations and
$\S_n(2\mn1\mn3,\ol{2}^o\mn31)$ can be established by composing our
bijection into lattice paths with the one from~\cite{BEM}.

\ms

Analogously to the definition for the pattern $\ol{2}^o\mn31$, we
say that a permutation $\pi$ avoids the pattern $\ol{2}^e\mn31$ if
every occurrence of $21$ in $\pi$ is part of an even number of
occurrences of $2\mn31$. We can also enumerate
$\{2\mn1\mn3,\ol{2}^e\mn31\}$-avoiding permutations.

\begin{prop}\label{th:even}
The number of $\{2\mn1\mn3,\ol{2}^e\mn31\}$-avoiding permutations of
size $n$ is
$$\frac{1}{n}\sum_{k=0}^{\lfloor n/2\rfloor}\left[2\binom{n}{2k}\binom{n-k}{k-1}+\frac{n}{n-k}\binom{n}{2k+1}\binom{n-k}{k}\right].$$
\end{prop}

\begin{proof}
Let $Q(t)=\sum_{n\ge1}|\S_n(2\mn1\mn3,\ol{2}^e\mn31)|\ t^n$. An
argument similar to the proof of Proposition~\ref{th:odd} shows that
$$Q(t)=\frac{(2-4t)\tilde{f}(t)^2+(-2+12t-7t^2-\tilde{g}(t))\tilde{f}(t)+(2-8t)\tilde{g}(t)+8t^3+46t^2-8t-4}{3t\tilde{f}(t)^2},$$
where $\tilde{g}(t)=\sqrt{3(-5t^4+24t^3-4t^2+12t-4)}$ and
$\tilde{f}(t)=[4(3t\tilde{g}(t)-11t^3-12t^2-6t+2)]^{1/3}$. It
follows that $Q=Q(t)$ is a root of the polynomial
$tQ^3+(4t-2)Q^2+(4t-1)Q+t=0$. Applying Lagrange inversion we get the
stated formula.
\end{proof}

\subsection{$\{2\mn1\mn3,2\mn3\mn41,3\mn2\mn41\}$-avoiding
permutations} The rightward generating tree for this class of
permutations has a simple succession rule. This allows us to
enumerate them easily. Let $K(t,u)=\sum_{n\ge1}\
\sum_{\pi\in\S_n(2\mn1\mn3,2\mn3\mn41,3\mn2\mn41)} u^{r(\pi)}
t^n=\sum_{r\ge1} K_r(t) u^r$.

\begin{prop}
The generating function for
$\{2\mn1\mn3,2\mn3\mn41,3\mn2\mn41\}$-avoiding permutations where
$u$ marks the value of the rightmost entry is
$$K(t,u)=\frac{1-t-2tu-\sqrt{1-2t-3t^2}}{2t(\frac{1}{u}+1+u)-2}.$$
\end{prop}

\begin{proof}
The succession rule for this class of permutations is \bce\bt{l}
$(r)\longrightarrow \begin{cases} (1)\ (2) & \mbox{if } r=1, \\
(r-1)\ (r)\ (r+1) & \mbox{if } r>1, \end{cases}$ \et\ece with the
root labeled $(1)$. This translates into the functional equation
\beq\label{eq:ker_K}\left[1-t\left(\frac{1}{u}+1+u\right)\right]K(t,u)=tu-tK_1(t).\eeq
Applying the kernel method we find that
$K_1(t)=\frac{1-t-\sqrt{1-2t-3t^2}}{2t}$, and substituting back into
(\ref{eq:ker_K}) we get the expression for $K(t,u)$.
\end{proof}

The generating function $K(t,1)$ also enumerates
$\{1\mn3\mn2,123\mn4\}$-avoiding permutations, as shown in
\cite[Example 2.6]{Man}. However, no direct bijection between
$\S_n(2\mn1\mn3,2\mn3\mn41,3\mn2\mn41)$ and
$\S_n(1\mn3\mn2,123\mn4)$ seems to be known.

\section{Generating trees with two labels}\label{sec:twolab}

The generating trees in all the examples in the previous section
were described using one label for each node. This will not be the
case in the families of permutations in this section. However, we
will use the same technique of translating the succession rule into
a set of functional equations and applying the kernel method to
them. This method is what unifies the different classes of
permutations studied in this paper.

Here we enumerate some classes of permutations whose rightward
generating tree has a succession rule that can be described using a
pair of labels for each node. These trees give rise to functional
equations with three variables. Even though no method is known to
solve them in general, in this section we present special cases
where we have been able to solve the corresponding equations.

In a few cases, one of the two labels is the length of the
permutation. When that happens, the functional equations have only
two variables, but the variable $t$ appears multiplied by another
variable, which makes these equations more difficult than the ones
in Section~\ref{sec:onelab}.

Note that for the classes that we consider in this section, the
enumeration of the permutations by their length has already been
done by different authors~\cite{C,ClaMan,EliMan,EliNoy,Man,Wes96}.
Our contribution is a refined enumeration of these permutations by
several parameters, and also the fact that our results are obtained
using the unifying framework of rightward generating trees.

\subsection{$\{2\mn1\mn3,12\mn3\}$-avoiding permutations}

It was shown in \cite{C} that $|\S_n(1\mn3\mn2,1\mn23)|=M_n$. A
bijection between $\S_n(1\mn3\mn2,1\mn23)$ and the set of Motzkin
paths of length $n$ was given in \cite{EliMan}. Clearly the sets
$\S_n(1\mn3\mn2,1\mn23)$ and $\S_n(2\mn1\mn3,12\mn3)$ are
equinumerous, since a permutation $\pi_1\pi_2\cdots\pi_n$ is
$\{1\mn3\mn2,1\mn23\}$-avoiding exactly when
$(n+1-\pi_n)\cdots(n+1-\pi_2)(n+1-\pi_1)$ is
$\{2\mn1\mn3,12\mn3\}$-avoiding. In this section we recover the
formula for $|\S_n(2\mn1\mn3,12\mn3)|$ using a generating tree with
two labels. This method provides a refined enumeration of
$\{2\mn1\mn3,12\mn3\}$-avoiding permutations by two new parameters:
the value of the last entry and the smallest value of the top of an
ascent.

Let $\T_1$ be the rightward generating tree for the set of
$\{2\mn1\mn3,12\mn3\}$-avoiding permutations. Given any
$\pi\in\S_n$, define the parameter \beq\label{eq:lab_lr}
\ell(\pi)=\begin{cases} n+1 & \mbox{if } \pi=n(n-1)\cdots21, \\
\min\{\pi_i : i>1,\ \pi_{i-1}<\pi_i\} & \mbox{otherwise.}
\end{cases} \eeq

Let each permutation $\pi$ be labeled by the pair
$(\ell,r)=(\ell(\pi),r(\pi))$. Note that since $\pi$ avoids
$12\mn3$, then necessarily $\ell\ge r$.

\begin{lemma}\label{lemma:gt_motzkin} The rightward generating tree $\T_1$ for $\{2\mn1\mn3,12\mn3\}$-avoiding
permutations is specified by the following succession rule on the
labels:
\bce\bt{l}
$(2,1)$ \\
$(\ell,r)\longrightarrow\begin{cases} (\ell+1,1)\ (\ell+1,2)\
\cdots\ (\ell+1,\ell) & \mbox{if } \ell=r,
\\ (\ell+1,1)\ (\ell+1,2)\ \cdots\ (\ell+1,r)\ (r+1,r+1) & \mbox{if } \ell>r.
\end{cases}$
\et\ece
\end{lemma}

\begin{proof} The permutation obtained by appending an entry to the right of
$\pi\in\S_n(2\mn1\mn3,12\mn3)$ is $2\mn1\mn3$-avoiding if and only
if the appended entry is at most $r(\pi)+1$, and it is
$12\mn3$-avoiding if and only if the appended entry is at most
$\ell(\pi)$. The labels of the children are obtained by looking at
how the values of $(\ell,r)$ change when the new entry is added.
\end{proof}

We will use this generating rule to obtain a formula for the
generating function
$$M(t,u,v):=\sum_{n\ge1}\ \sum_{\pi\in\S_n(2\mn1\mn3,12\mn3)}
u^{\ell(\pi)} v^{r(\pi)}\ t^n.$$ For fixed $\ell$ and $r$, let
$M_{\ell,r}(t)=\sum_{n\ge1}
|\{\pi\in\S_n(2\mn1\mn3,12\mn3):\ell(\pi)=\ell,r(\pi)=r\}|\ t^n$.
Note that $M(t,u,v)=\sum_{\ell,r} M_{\ell,r}(t) u^\ell v^r$.

\begin{prop}
The generating function for $\{2\mn1\mn3,12\mn3\}$-avoiding
permutations where $u$ and $v$ mark the parameters $\ell$ and $r$
defined above is
$$M(t,u,v)=\frac{[(1-u)v+c_1t+c_2t^2+c_3t^3+c_4t^4-((1-u)v+tu+t^2u^2v)\sqrt{1-2t-3t^2})]u^2v}{2(1-u-tu(1-u)+t^2u^2)(1-uv+tuv+t^2u^2v^2)},$$
where $c_1=2-u-v-uv+2u^2v$, $c_2=u(-1+(2-u)v+2(u-1)v^2)$,
$c_3=u^2v(-3+2v-2uv)$, and $c_4=-2u^3v^2$.
\end{prop}

Substituting $u=v=1$ in the above expression we recover the
generating function for the Motzkin numbers.

\begin{proof}
The coefficient of $t^n$ in $M(t,u,v)$ is the sum of $u^\ell v^r$
over all the pairs $(\ell,r)$ of labels that appear at level $n$ of
the tree. By Lemma~\ref{lemma:gt_motzkin}, the children of a node
with labels $(\ell,\ell)$ contribute
$u^{\ell+1}v+u^{\ell+1}v^2+\cdots+u^{\ell+1}v^\ell$ to the next
level, and the children of a node with labels $(\ell,r)$ with
$\ell>r$ contribute
$u^{\ell+1}v+u^{\ell+1}v^2+\cdots+u^{\ell+1}v^r+u^{r+1}v^{r+1}$. It
follows that \beq\label{eq:Mtuv} M(t,u,v) = tu^2v+t \sum_\ell
M_{\ell,\ell}(t)u^{\ell+1}(v+v^2+\cdots+v^\ell) + t \sum_{\ell>r}
M_{\ell,r}(t)[u^{\ell+1}(v+v^2+\cdots+v^r)+u^{r+1}v^{r+1}].\eeq It
will be convenient to define
$$M_>(t,u,v):=\sum_{n\ge1}\ \underset{\mathrm{with\ }\ell(\pi)>r(\pi)}{\sum_{\pi\in\S_n(2\mn1\mn3,12\mn3)}}
u^{\ell(\pi)} v^{r(\pi)}\ t^n \quad\mbox{and}\quad
M_=(t,u,v):=\sum_{n\ge1}\ \underset{\mathrm{with\
}\ell(\pi)=r(\pi)}{\sum_{\pi\in\S_n(2\mn1\mn3,12\mn3)}}
(uv)^{\ell(\pi)}\ t^n,$$ so that $M(t,u,v)=M_>(t,u,v)+M_=(t,u,v)$.
Taking from (\ref{eq:Mtuv}) only the pairs $(\ell,r)$ with $\ell>r$,
we get \bea\nn M_>(t,u,v)&=&tu^2v+t \sum_\ell
M_{\ell,\ell}(t)u^{\ell+1}(v+v^2+\cdots+v^\ell)+t
\sum_{\ell>r} M_{\ell,r}(t)[u^{\ell+1}(v+v^2+\cdots+v^r)] \\
&=& tu^2v+t \sum_\ell
M_{\ell,\ell}(t)u^{\ell+1}\frac{v^{\ell+1}-v}{v-1}+t
\sum_{\ell>r} M_{\ell,r}(t)u^{\ell+1}\frac{v^{r+1}-v}{v-1}\nn \\
&=&
tu^2v+\frac{tuv}{v-1}\left[M_=(t,u,v)-M_=(t,u,1)+M_>(t,u,v)-M_>(t,u,1)\right].\label{eq:Mgtuv}\eea
Similarly, taking from (\ref{eq:Mtuv}) only the pairs $(\ell,r)$
with $\ell=r$, \beq M_=(t,u,v)=t \sum_{\ell>r}
M_{\ell,r}(t)u^{r+1}v^{r+1}=tuv\sum_{\ell>r}
M_{\ell,r}(t)(uv)^r=tuv\ M_>(t,1,uv).\label{eq:Metuv}\eeq

Using in (\ref{eq:Mgtuv}) the expression of $M_=$ in terms of $M_>$
given in (\ref{eq:Metuv}), we get \beq
M_>(t,u,v)=tu^2v+\frac{tuv}{v-1}\left[tuv\ M_>(t,1,uv)-tu\
M_>(t,1,u)+M_>(t,u,v)-M_>(t,u,1)\right].\label{eq:Mg}\eeq
Substituting $u=1$ in this equation and collecting the terms in
$M_>(t,1,v)$, we have
\beq\left(1-\frac{t^2v^2}{v-1}-\frac{tv}{v-1}\right)M_>(t,1,v)=tv-\frac{t(t+1)v}{v-1}
M_>(t,1,1).\label{eq:kernel_Mgt1v}\eeq

Now we apply the kernel method, substituting
$v=v_0=v_0(t)=\frac{1-t-\sqrt{1-2t-3t^2}}{2t^2}$ in
(\ref{eq:kernel_Mgt1v}) to obtain
$$M_>(t,1,1)=\frac{v_0-1}{t+1}=\frac{1-t-2t^2-\sqrt{1-2t-3t^2}}{2t^2(t+1)}.$$
Plugging this expression for $M_>(t,1,1)$ back into
(\ref{eq:kernel_Mgt1v}) we get that \beq\label{eq:Mgt1v}
M_>(t,1,v)=\frac{(1-t-2t^2v-\sqrt{1-2t-3t^2})v}{2t(1-v+tv+t^2v^2)}.\eeq

If we write equation (\ref{eq:Mg}) as
\beq\label{eq:kernel_Mgtuv}\left(1-\frac{tuv}{v-1}\right)M_>(t,u,v)=tu^2v+\frac{tuv}{v-1}\left[tuv\
M_>(t,1,uv)-tu\ M_>(t,1,u)-M_>(t,u,1)\right],\eeq we can apply again
the kernel method, taking $v=v_1=v_1(t,u)=\frac{1}{1-tu}$.
This cancels the left hand side and gives \beq\nn M_>(t,u,1)=
\frac{[2(1-u)+u^2-t(1+2t)u^2+(1-2u)\sqrt{1-2t-3t^2})]tu^2}{2(1-u+tu+t^2u^2)(1-u-tu(1-u)+t^2u^2)}\eeq
using (\ref{eq:Mgt1v}). Substituting back into
(\ref{eq:kernel_Mgtuv}) and using (\ref{eq:Mgt1v}) again we get that
$$M_>(t,u,v)=\frac{[2-u-uv+u^2v+tu(v-1)-t(1+2t)u^2v+(1-2u)\sqrt{1-2t-3t^2})]tu^2v}{2(1-u-tu(1-u)+t^2u^2)(1-uv+tuv+t^2u^2v^2)}.$$
Finally, combining it with the fact that
$$M(t,u,v)=M_>(t,u,v)+M_=(t,u,v)=M_>(t,u,v)+tuv\ M_>(t,1,uv),$$ we
obtain the desired expression for $M(t,u,v)$.
\end{proof}

We have encountered two classes of pattern-avoiding permutations
enumerated by the Motzkin numbers, namely \beq\label{eq:added}
|\S_{n+1}(2\mn1\mn3,\ol{2}\mn31)|=|\S_n(2\mn1\mn3,12\mn3)|=M_n \eeq
(see Proposition~\ref{th:213bar231}). In Section~\ref{sec:213bar231}
we described a bijection $\varphi$ between $2\mn3\mn1$-avoiding
permutations and Dyck paths. A permutation $\pi$ is
$\ol{2}\mn31$-avoiding if and only if $\varphi(\pi)$ is a $UDU$-free
Dyck path. It is not hard to check (see~\cite{EliMan}) that $\pi$ is
$12\mn3$-avoiding if and only if $\varphi(\pi)$ is $UUU$-free. Next
we give a bijection between $UDU$-free Dyck paths of size $n+1$ and
$UUU$-free Dyck paths of size $n$, reproving
equation~(\ref{eq:added}).

Given a $UDU$-free Dyck path, mark each $D$ that is immediately
preceded and followed by a $D$ (and also the rightmost $D$ if it is
preceded by a $D$). Move left each one of the marked $D$s so that it
immediately follows its matching $U$. Finally, delete the rightmost
peak (i.e., occurrence of $UD$). This gives a $DDD$-free Dyck path,
which can be easily turned into a $UUU$-free one by reversing it,
that is, reading the steps from right to left and exchanging $U$s
and $D$s.

To show that this map is the desired bijection, we now describe its
inverse. Given a $DDD$-free Dyck path, we first reverse it and then
append a peak $UD$ to it. Define the height of a step to be the
number of $U$s minus the number of $D$s preceding it. Mark each $D$
step in an occurrence of $UDU$. For each marked step, if $h$ is its
height, move it to the right so that it immediately precedes the
next $D$ step with height $h-1$ (if $h=1$, then move it to the end).
This produces the original $UDU$-free Dyck path.

As an example of this bijection, consider the $UDU$-free path
$U\bar{U}UUDD\bar{U}UD\bar{D}\bar{D}D\bar{U}UD\bar{D}$. The marked
$D$s and their matching $U$s are distinguished with a bar. The
$DDD$-free Dyck path that we obtain is
$UU\bar{D}UUDDU\bar{D}UDDU\bar{D}$ (the barred $D$s are the steps
that have been moved), and its reversal is $UDUUDUDUUDDUDD$.
Applying the inverse map moves the barred $D$s back to their
original position.

\subsection{$\{2\mn1\mn3,32\mn1\}$-avoiding
permutations}

It is known \cite{C} that $|\S_n(2\mn1\mn3,32\mn1)|=2^{n-1}$. Here
we use rightward generating trees with two labels to recover this
fact, and to refine it with two parameters: the value of the last
entry and the largest value of the bottom of a descent. Given any
$\pi\in\S_n$, define \beq\nn h(\pi)=\begin{cases} 0 & \mbox{if }
\pi=12\cdots n,
\\ \max\{\pi_i : i>1,\ \pi_{i-1}>\pi_i\} & \mbox{otherwise.}
\end{cases} \eeq

To each $\{2\mn1\mn3,32\mn1\}$-avoiding permutation $\pi$ we assign
the pair of labels $(h,r)=(h(\pi),r(\pi))$. Note that since $\pi$
avoids $32\mn1$, then necessarily $h\le r$.

\begin{lemma}\label{lemma:gt_2} The rightward generating tree for $\{2\mn1\mn3,32\mn1\}$-avoiding
permutations is specified by the following succession rule on the
labels: \bce\bt{l}
$(0,1)$ \\
$(h,r)\longrightarrow (h+1,h+1)\ (h+2,h+2)\ \cdots\ (r,r)\ (h,r+1).$
\et\ece
\end{lemma}

\begin{proof} When we append an entry $i$ to a $\{2\mn1\mn3,32\mn1\}$-avoiding
permutation, the new permutation is $2\mn1\mn3$-avoiding if and only
if $i\le r(\pi)+1$, and it is $32\mn1$-avoiding if and only if $i>h(\pi)$. The list of labels of the children obtained by appending an
$i$ satisfying these two conditions is the right hand side of the
rule.
\end{proof}

Let
$$N(t,u,v)=\sum_{n\ge1}\ \sum_{\pi\in\S_n(2\mn1\mn3,32\mn1)}
u^{h(\pi)} v^{r(\pi)}\ t^n=\sum_{h,r} N_{h,r}(t) u^h v^r.$$

\begin{prop}
The generating function for $\{2\mn1\mn3,32\mn1\}$-avoiding
permutations where $u$ and $v$ mark the parameters $h$ and $r$
defined above is
$$N(t,u,v)=\frac{tv(1-t+tu-tuv)}{(1-tv)(1-t-tuv)}.$$
\end{prop}

\begin{proof}
By Lemma~\ref{lemma:gt_2}, the children of a node with labels
$(h,r)$ contribute $u^{h+1}v^{h+1}+u^{h+2}v^{h+2}+\cdots+u^r v^r+u^h
v^{r+1}$ to the next level. It follows that \bea\nn N(t,u,v)&=&tv +
t \sum_{h,r} N_{h,r}(t)\left[\frac{(uv)^{r+1}-(uv)^{h+1}}{uv-1}+u^h
v^{r+1}\right]\\ &=& tv + tv N(t,u,v)+\frac{tuv
[N(t,1,uv)-N(t,uv,1)]}{uv-1}.\label{eq:Ntuv}\eea

Substituting $u=1$ and $v=1$ separately gives a system of two
equations in $N(t,1,*)$ and $N(t,*,1)$ that can be easily solved.
\end{proof}

The above result can indeed be obtained as well without using
rightward generating trees. The recursive structure of
$2\mn1\mn3$-avoiding permutations (i.e., they are of the form
$\sigma1\tau$, where $\sigma$ and $\tau$ are $2\mn1\mn3$-avoiding
and every entry in $\sigma$ is larger than every entry in $\tau$) can be
used to obtain an equation satisfied by $N(t,u,v)$ and to deduce the
above formula without much difficulty.

\subsection{$\{2\mn1\mn3,34\mn21\}$-avoiding
permutations}

The labels that will be convenient to use to describe the rightward
generating tree for this class are $(s,r)=(s(\pi),r(\pi))$, where
\beq\label{eq:defs} s(\pi)=\begin{cases} 0 & \mbox{if }
\pi=n(n-1)\cdots 21,
\\ \max\{\pi_i : \ \pi_{i}<\pi_{i+1}\} & \mbox{otherwise,}
\end{cases} \eeq
and $r(\pi)=\pi_n$ as usual.

\begin{lemma}\label{lemma:3421} The rightward generating tree for $\{2\mn1\mn3,34\mn21\}$-avoiding
permutations is specified by the following succession rule on the
labels: \bce\bt{l}
$(0,1)$ \\
$(s,r)\longrightarrow\begin{cases} (s+1,1)\ (s+1,2)\ \cdots\
(s+1,s)\ (s,s+1)\ (r,r+1) & \mbox{if } s<r,
\\ (s+1,r+1) & \mbox{if } s>r.
\end{cases}$
\et\ece
\end{lemma}

\begin{proof}
First note that the $2\mn1\mn3$-avoiding condition implies that the new
entry appended to $\pi$ has to be at most $r+1$. If $s<r$ and
$\pi_i$ is an entry to the right of $s$, then $\pi_i>s$, otherwise
$s\pi_i r$ would be an occurrence of $2\mn1\mn3$. In fact, we also
know that $\pi_i\neq s+1$, unless $\pi_i=r=s+1$, because otherwise
the entry following $s+1$ would be greater than it, contradicting
the definition of $s$. So, unless $r=s+1$, the entry $s+1$ precedes
$s$, so the appended entry cannot be greater than $s+1$, otherwise
it would create a $2\mn1\mn3$. This explains the labels in the case
$s<r$. If $s>r$, the appended entry has to be greater than $r$ for
the new permutation to be $34\mn21$-avoiding.
\end{proof}

Let
$$K(t,u,v):=\sum_{n\ge1}\ \sum_{\pi\in\S_n(2\mn1\mn3,34\mn21)}
u^{s(\pi)} v^{r(\pi)}\ t^n=\sum_{s,r} K_{s,r}(t) u^s v^r,$$ and let
$K_<(t,u,v)$ and $K_>(t,u,v)$ be defined similarly, with the sum
running only over permutations with $s(\pi)<r(\pi)$ and
$s(\pi)>r(\pi)$, respectively, so that
$K(t,u,v)=K_<(t,u,v)+K_>(t,u,v)$.

\begin{prop}
The generating function for $\{2\mn1\mn3,34\mn21\}$-avoiding
permutations where $u$ and $v$ mark the parameters $s$ and $r$
defined above is
\beq\label{eq:K}K(t,u,v)=\frac{tv[1-(1+u+uv)t+(u^2+uv+u^2v)t^2]}{(1-t-tu)(1-t-tuv)(1-tuv)}.\eeq
\end{prop}

\begin{proof}
By Lemma~\ref{lemma:3421}, the generating functions $K_<$ and $K_>$
satisfy \bea\nn K_<(t,u,v)&=& tv+t\sum_{s<r}
K_{s,r}(t)(u^sv^{s+1}+u^rv^{r+1})= tv+tv[K_<(t,uv,1)+K_<(t,1,uv)],
\\ \nn
K_>(t,u,v)&=&
t\sum_{s<r}K_{s,r}(t)u^{s+1}(v+\cdots+v^s)+t\sum_{s>r}K_{s,r}(t)u^{s+1}v^{r+1}\\
\nn &=& \frac{tuv}{v-1}[K_<(t,uv,1)-K_<(t,u,1)]+tuv K_>(t,u,v).\eea

Substituting first $u=1$ and then $v=1$ in the first equation, we
get two equations involving $K_<(t,1,w)$ and $K_<(t,w,1)$ that can
be easily solved to give
$$K_<(t,u,v)=\frac{tv}{1-t-tuv}.$$ The second equation then implies
that
$$K_>(t,u,v)=\frac{u^2vt^3}{(1-t-tu)(1-t-tuv)(1-tuv)},$$ and
the proposition follows.
\end{proof}

\begin{corollary}\label{cor:west}
The number of $\{2\mn1\mn3,34\mn21\}$-avoiding permutations of size
$n$ is $(n-1)2^{n-2}+1$.
\end{corollary}

\begin{proof}
Taking $u=v=1$ in (\ref{eq:K}), we get that
$$K(t,1,1)=\frac{t(1-3t+3t^2)}{(1-t)(1-2t)^2}.$$
The coefficient of $t^n$ in the series expansion of this rational
function is $(n-1)2^{n-2}+1$.
\end{proof}

It is not hard to show that
$\S_n(2\mn1\mn3,34\mn21)=\S_n(2\mn1\mn3,3\mn4\mn2\mn1)=\S_n(1\mn3\mn2,3\mn4\mn2\mn1)$.
This last set of permutations was enumerated by West~\cite{Wes96},
and Corollary~\ref{cor:west} agrees with his result.

\subsection{$\{1\mn2\mn34,2\mn1\mn3\}$-avoiding permutations}

The generating function for these permutations appears
in~\cite{Man}. In fact, it is easy to see that
$\S_n(1\mn2\mn34,2\mn1\mn3)=\S_n(1\mn2\mn3\mn4,2\mn1\mn3)$, and the
latter set of permutations was studied in~\cite{Wes96}, where it is
shown that they are counted by the Fibonacci numbers $F_{2n-1}$.
Here we derive the generating function and obtain a refinement of it
using a rightward generating tree with two labels.

Let the labels of a permutation $\pi$ be the pair
$(m,r)=(m(\pi),r(\pi))$, where $r(\pi)=\pi_n$ and
$$ m(\pi)=\begin{cases} n+1 & \mbox{if } \pi=n(n-1)\cdots21, \\ \min\{\pi_i : \exists j<i \mbox{ with } \pi_{j}<\pi_i\} & \mbox{otherwise.}
\end{cases} $$
Note that we always have $m(\pi)\le r(\pi)$ unless $r=1$, and that
if $m(\pi)=r(\pi)$, then $\pi=n(n-1)\cdots312$, so $m=r=2$.

\begin{lemma}\label{lemma:gt_mansour2} The rightward generating tree for $\{1\mn2\mn34,2\mn1\mn3\}$-avoiding
permutations is specified by the following succession rule on the
labels: \bce\bt{l}
$(2,1)$ \\
$(m,r)\longrightarrow\begin{cases} (m+1,1)\ (2,2) & \mbox{if } r=1,
\\ (3,1)\ (2,2)\ (2,3)
& \mbox{if } m=r=2,
\\ (m+1,1)\ (2,2) \ (m,m+1)\ \cdots\ (m,r) & \mbox{if }
m<r.
\end{cases}$
\et\ece
\end{lemma}

\begin{proof}
As usual, the appended entry has to be at most $r+1$ for the
permutation to avoid $2\mn1\mn3$. In the case that $m<r$, this entry
cannot be greater than $r$ in order to avoid $1\mn2\mn34$. The
labels are now obtained by looking at how the parameter $m$ changes
after appending the new entry.
\end{proof}

Let $H(t,u,v):=\sum_{n\ge1}\ \sum_{\pi\in\S_n(1\mn2\mn34,2\mn1\mn3)}
u^{m(\pi)} v^{r(\pi)}\ t^n$, and let $H_1(t,u,v)$, $H_=(t,u,v)$, and
$H_<(t,u,v)$ be defined similarly, with the summation restricted to
permutations with $r(\pi)=1$, $m(\pi)=r(\pi)$, and $m(\pi)<r(\pi)$,
respectively, so that $H(t,u,v)=H_1(t,u,v)+H_=(t,u,v)+H_<(t,u,v)$.

\begin{prop}
The generating function for $\{1\mn2\mn34,2\mn1\mn3\}$-avoiding
permutations where $u$ and $v$ mark the parameters $m$ and $r$
defined above is \beq \nn
H(t,u,v)=\frac{tu^2v[1+(v-3)t+(1+u-v-uv+v^2)t^2+uv(1-v)t^3]}{(1-3t+t^2)(1-tu)}.\eeq
\end{prop}

\begin{proof}
From Lemma~\ref{lemma:gt_mansour2} we get the following functional
equations defining $H_1$, $H_=$, and $H_<$. \bea \label{eqH1}
H_1(t,u,v)&=&tu^2v+tuv H(t,u,1)
\\ \label{eqHe}
H_=(t,u,v)&=&tu^2v^2 H(t,1,1)
\\ \label{eqHl}
H_<(t,u,v)&=&tv
H_=(t,u,v)+\frac{tv}{v-1}[H_<(t,u,v)-H_<(t,uv,1)]\eea
 Combining (\ref{eqHe}) and (\ref{eqHl}), introducing a
variable $w=uv$, and defining $\H_<(t,w,v)=H_<(t,\frac{w}{v},v)$, we
get \beq\label{eqHlH} \left(1-\frac{tv}{v-1}\right)\H_<(t,w,v)=t^2 v
w^2 H(t,1,1)-\frac{tv}{v-1}\H_<(t,w,1).\eeq
 The kernel is canceled with $v=\frac{1}{1-t}$, giving an expression
for $\H_<(t,w,1)$ in terms of $H(t,1,1)$, which plugged back into
(\ref{eqHlH}) yields
\beq\label{eqHlH2}H_<(t,u,v)=\H_<(t,uv,v)=\frac{t^2u^2v^3}{1-t}H(t,1,1).\eeq

On the other hand, adding equations (\ref{eqH1}) and (\ref{eqHe})
and using that $H_1(t,u,v)+H_=(t,u,v)=H(t,u,v)-H_<(t,u,v)$, we get
$$H_<(t,u,v) = H(t,u,v)-tu^2v-tuvH(t,u,1)-tu^2v^2H(t,1,1),$$
which combined with (\ref{eqHlH2}) gives a simple expression
relating $H(t,u,v)$, $H(t,u,1)$ and $H(t,1,1)$. In this expression,
the substitution $u=v=1$ gives
$$H(t,1,1)=\frac{t(1-t)}{1-3t+t^2},$$ and the substitution $v=1$
puts $H(t,u,1)$ in terms of $H(t,1,1)$. All together produces the
desired formula for $H(t,u,v)$.
\end{proof}

\subsection{$\{12\mn34,2\mn1\mn3\}$-avoiding permutations}

It was proved in \cite{Man} that the generating function for
permutations avoiding $\{12\mn34,2\mn1\mn3\}$ is
$\frac{1-2t-t^2-\sqrt{1-4t+2t^2+t^4}}{2t^2}$. Using the labels
$(\ell,r)$ defined as in (\ref{eq:lab_lr}), we can construct a
generating tree with two labels for this class of permutations. The
proof of the following lemma is straightforward and analogous to
that of Lemma~\ref{lemma:gt_motzkin}.

\begin{lemma}\label{lemma:gt_mansour} The rightward generating tree for $\{12\mn34,2\mn1\mn3\}$-avoiding
permutations is specified by the following succession rule on the
labels: \bce\bt{l}
$(2,1)$ \\
$(\ell,r)\longrightarrow\begin{cases} (\ell+1,1)\ (\ell+1,2)\
\cdots\ (\ell+1,r)\ (r+1,r+1) & \mbox{if } \ell>r,
\\ (\ell+1,1)\ (\ell+1,2)\ \cdots\ (\ell+1,\ell)\ (\ell,\ell+1)
& \mbox{if } \ell=r,
\\ (\ell+1,1)\ (\ell+1,2)\ \cdots\ (\ell+1,\ell)\ (\ell,\ell+1)\ (\ell,\ell+2)\ \cdots\ (\ell,r) & \mbox{if }
\ell<r.
\end{cases}$
\et\ece
\end{lemma}

Let $F(t,u,v):=\sum_{n\ge1}\ \sum_{\pi\in\S_n(12\mn34,2\mn1\mn3)}
u^{\ell(\pi)} v^{r(\pi)}\ t^n$, and let $F_>(t,u,v)$, $F_=(t,u,v)$,
and $F_<(t,u,v)$ be defined similarly, with the summation restricted
to permutations with $\ell(\pi)>r(\pi)$, $\ell(\pi)=r(\pi)$, and
$\ell(\pi)<r(\pi)$, respectively. By definition,
$F(t,u,v)=F_>(t,u,v)+F_=(t,u,v)+F_<(t,u,v)$.

\begin{prop}
The generating function for $\{12\mn34,2\mn1\mn3\}$-avoiding
permutations where $u$ and $v$ mark the parameters $\ell$ and $r$
defined above is \beq \nn
F(t,u,v)=\frac{u^2v[p_1(t,u,v)+p_2(t,u,v)\sqrt{1-4t+2t^2+t^4}]}{2[(1+tuv)^2-uv-t-uvt^2][1+(u+t)(tu-1)]},\eeq
where
 {\small
\bea\nn
p_1(t,u,v)&=&(1-u)v+(2-u-4v+2uv+v^2+2u^2v-uv^2)t+(-4+u+6v+uv-3v^2-6u^2v+3u^2v^2)t^2\\
\nn &&+\,(2+u-4v-5uv+3v^2+4u^2v+4uv^2-4u^2v^2-2uv^3-2u^3v^2+2u^2v^3)t^3\\
\nn
&&+\,(-u+v+4uv-v^2-4uv^2-u^2v^2+2uv^3+2u^3v^2-2u^3v^3)t^4-uv(v-1)(2uv-1)t^5,\\
\nn
p_2(t,u,v)&=&(u-1)v+[(u-1)v(v-2)-u]t+(u-v+v^2-u^2v^2)t^2+uv(1-v)t^3.\eea}
\end{prop}

Note that this expression becomes much simpler if we ignore the
parameter $r$, that is,
$$F(t,u,1)=\frac{u^2(1-2tu-t^2-\sqrt{1-4t+2t^2+t^4})}{2[1+(u+t)(tu-1)]},$$
and coincides with the result from~\cite{Man} if we ignore both
parameters:
$$F(t,1,1)=\frac{1-2t-t^2-\sqrt{1-4t+2t^2+t^4}}{2t^2}.$$

\begin{proof}
From Lemma~\ref{lemma:gt_mansour} we get the following functional
equations defining $F_>$, $F_=$, and $F_<$. \bea\label{eqF>}
F_>(t,u,v)&=&tu^2v+\frac{tuv}{v-1}[F_>(t,u,v)-F_>(t,u,1)\\
\nn &&\hspace{21mm}+\ F_=(t,u,v)-F_=(t,u,1)+F_<(t,uv,1)-F_<(t,u,1)],
\\ \label{eqF=}
F_=(t,u,v)&=&tuv F_>(t,1,uv),
\\ \label{eqF<}
F_<(t,u,v)&=&tv
F_=(t,u,v)+\frac{tv}{v-1}[F_<(t,u,v)-F_<(t,uv,1)].\eea

We can introduce a variable $w=uv$ in (\ref{eqF<}), and apply the
kernel method with $v=\frac{1}{1-t}$ to get that
$$F_<(t,w,1)=\frac{t^2 w F_>(t,1,w)}{1-t}.$$
Using this expression together with (\ref{eqF=}) in (\ref{eqF>}), we
get an equation that involves only $F_>$:
\beq\label{eqF>uv}\left(1-\frac{tuv}{v-1}\right)F_>(t,u,v)=tu^2v-\frac{tuv}{v-1}\left[F_>(t,u,1)+\frac{tuv}{1-t}F_>(t,1,uv)-
\frac{tu}{1-t}F_>(t,1,u)\right].\eeq Substituting $u=1$, it becomes
\beq\label{eqF>1v}\frac{(1-v-t+2tv-vt^2+v^2t^2)}{(1-v)(1-t)}F_>(t,1,v)=tv+\frac{tv}{(1-v)(1-t)}F_>(t,1,1).\eeq
We apply the kernel method again, this time with
$v=\frac{1-2t+t^2-\sqrt{1-4t+2t^2+t^4}}{2t^2}$ to cancel the left
hand side of (\ref{eqF>1v}), which yields
$$F_>(t,1,1)=\frac{1-3t+t^2+t^3+(t-1)\sqrt{1-4t+2t^2+t^4}}{2t^2}.$$
Now we can use (\ref{eqF>1v}) to obtain a formula for $F_>(t,1,v)$.
Applying again the kernel method in (\ref{eqF>uv}), with
$v=\frac{1}{1-tu}$, we get an expression for $F_>(t,u,1)$ in terms
of $F_>(t,1,u)$ and $F_>(t,1,\frac{u}{1-tu})$, and therefore a
formula for $F_>(t,u,1)$. Substituting back into (\ref{eqF>uv}), we
get a formula for $F_>(t,u,v)$. From this it is straightforward to
obtain formulas for $F_<(t,u,v)$ and $F_=(t,u,v)$ as well, and the
result follows.
\end{proof}

\subsection{$\{1\mn23,3\mn12\}$-avoiding
permutations}\label{sec:1_23_3_12}

Generating trees with two labels can be used to obtain the
generating function for the number of $\{1\mn23,3\mn12\}$-avoiding
permutations. These permutations were studied in \cite{ClaMan},
where it was shown that if we let $b_n=|\S_n(1\mn23,3\mn12)|$, then
these numbers satisfy the recurrence
$b_{n+2}=b_{n+1}+\sum_{k=0}^n\binom{n}{k}b_k$. Here we obtain an
ordinary generating function without going through the recurrence.
The labels are particularly easy in this case because we can take
one of them to be just the length $n$ of the permutation. The labels
of $\pi\in\S_n$ are then $(r,n)$, where $r=\pi_n$ as usual. The
advantage of having one of the labels be $n$ is that we do not need
an extra variable for this label in the generating function, since
it is already encoded in the exponent of the variable $t$.

\begin{lemma}\label{lemma:gt_1_23_3_12} The rightward generating tree for $\{1\mn23,3\mn12\}$-avoiding
permutations is specified by the following succession rule on the
labels: \bce\bt{l}
$(1,1)$ \\
$(r,n)\longrightarrow\begin{cases} (1,n+1)\ (n+1,n+1) & \mbox{if }
r=1,
\\ (1,n+1)\ (2,n+1)\ \cdots\
(r,n+1) & \mbox{if } r>1.
\end{cases}$ \et\ece
\end{lemma}

\begin{proof}
The appended element cannot be larger than the rightmost entry of
$\pi$, except where this entry is $1$, in which case the appended
element can be the new largest one.
\end{proof}

Let $P(t,u):=\sum_{n\ge1}\ \sum_{\pi\in\S_n(1\mn23,3\mn12)}
u^{r(\pi)} \ t^n=\sum_{r} P_r(t) u^r$.

\begin{prop}
The generating function for $\{1\mn23,3\mn12\}$-avoiding
permutations is
$$P(t,1)=\sum_{k\ge1}\frac{t^{2k-1}(1-(k-1)t)}{(1-t)^2(1-2t)^2\cdots(1-kt)^2}.$$
\end{prop}

\begin{proof}
From Lemma~\ref{lemma:gt_1_23_3_12} we get \bea P(t,u) \label{eq:P}
= tu+\frac{tu}{u-1}(P(t,u)-u P_1(t)-P(t,1)+P_1(t))+tu
(P_1(t)+P_1(tu)).\eea Using that $P_1(t)=t+t\ P(t,1)$ and collecting
the terms with $P(t,u)$, we get
\beq\label{eq:Ptu}\left(1-\frac{tu}{u-1}\right)P(t,u)=tu+t^2u^2+t^2u^2
P(tu,1)+\left(t^2u+\frac{tu(t-1-tu)}{u-1}\right)P(t,1).\eeq

Substituting $u=\frac{1}{1-t}$ gives
$$P(t,1)=\frac{t}{(1-t)^2}\left(1+tP(\frac{t}{1-t},1)\right),$$
and by iterated application of this formula, \bea \nn
P(t,1)=\frac{t}{(1-t)^2}\left(1+\frac{t^2(1-t)}{(1-2t)^2}\left(1+\frac{t^2(1-2t)}{(1-t)(1-3t)^2}
\left(1+\frac{t^2(1-3t)}{(1-2t)(1-4t)^2}\left(1+\cdots\right)\right)\right)\right)\hspace{15mm}\\
\nn
=\frac{t}{(1-t)^2}+\frac{t^3}{(1-t)(1-2t)^2}+\frac{t^5}{(1-t)^2(1-2t)(1-3t)^2}+\frac{t^7}{(1-t)^2(1-2t)^2(1-3t)(1-4t)^2}+\cdots,\eea
which is the formula above. If we substitute this expression back
into (\ref{eq:Ptu}) we get the refined formula for $P(t,u)$.
\end{proof}

A very similar argument can be applied to $1\mn23$-avoiding
permutations, which are known to be enumerated by the Bell
numbers~\cite{C}. Our approach in this case gives essentially the
same functional equation that is derived in~\cite{FirMan} using what
the authors call the {\it scanning-elements algorithm}.

Rightward generating trees and the kernel method can also be used to
produce a functional equation for the ordinary generating function
of $123$-avoiding permutations. We omit this result here because a
more direct way to enumerate these permutations was already given
in~\cite{EliNoy}.

\section{Generating trees with three labels}\label{sec:threelab}

In this section we include two instances of permutations avoiding
generalized patterns where the rightward generating tree can be
described by a succession rule with three labels. One of these
labels is the length of the permutation, so that the functional
equations that we obtain have three variables instead of four.
However, the fact that the variable $t$ appears multiplied by
another variable adds some difficulty to the equations.

To the best of our knowledge, the two classes of restricted
permutations considered in this section have never been enumerated
before.

\subsection{$\{1\mn23,3\mn12,34\mn21\}$-avoiding permutations}
\label{sec:threelabsub1}

Given a permutation $\pi\in\S_n$, let $s(\pi)$ be defined as in
(\ref{eq:defs}). We associate to $\pi$ the triple of labels
$(s,r,n)=(s(\pi),r(\pi),n)$.

\begin{lemma}\label{lemma:1233123421} The rightward generating tree for $\{1\mn23,3\mn12,34\mn21\}$-avoiding
permutations is specified by the following succession rule on the
labels: \bce\bt{l}
$(0,1,1)$ \\
$(s,r,n)\longrightarrow\begin{cases} (s+1,1,n+1)\ (s+1,2,n+1)\
\cdots\ (s+1,s,n+1) & \\ \hspace{2cm} (s,s+1,n+1)\ (s,s+2,n+1)\
\cdots\ (s,r,n+1) \quad & \mbox{if } s<r\neq1,
\\ (0,1,n+1)\ (1,n+1,n+1) & \mbox{if } (s,r)=(0,1),
\\ (s,n+1,n+1) & \mbox{if } s>r=1,
\\ \emptyset & \mbox{if } s>r>1.
\end{cases}$
\et\ece
\end{lemma}

\begin{proof}
If $r>1$, the appended entry has to be at most $r$ for the new
permutation to avoid $1\mn23$. If $s>r$, it has to be at least $r+1$
for the new permutation to avoid $34\mn21$. Finally, if $r=1$, the
appended entry has to be $n+1$ for the permutation to avoid
$3\mn12$, unless $s=0$, which means that $\pi$ is the decreasing
permutation. Combining these conditions we get the four possible
cases and the new labels in each case.
\end{proof}

The four cases in the succession rule above suggest dividing the set
$\Theta$ of values that the pair $(s,r)$ can take into four disjoint
sets: $\Theta_1=\{(s,r):s<r\neq1\}$, $\Theta_2=\{(0,1)\}$,
$\Theta_3=\{(s,r):s>r=1\}$, $\Theta_4=\{(s,r):s>r>1\}$. For
$i=1,2,3,4$, let $$R_i(t,u,v):=\sum_{n\ge1}\ \underset{\mathrm{with\
}(s(\pi),r(\pi))\in\Theta_i}{\sum_{\pi\in\S_n(1\mn23,3\mn12,34\mn21)}}
u^{s(\pi)} v^{r(\pi)} \ t^n,$$ and let
$R(t,u,v)=R_1(t,u,v)+R_2(t,u,v)+R_3(t,u,v)+R_4(t,u,v)$.

\begin{prop}
The generating function for $\{1\mn23,3\mn12,34\mn21\}$-avoiding
permutations where $u$ marks the parameter $s$ defined above is
\beq\label{eq:R}1+R(t,u,1)=\sum_{k\ge0}\dfrac{t^{2k}u^k(1+ktu)}{(1-(k+1)t)\prod_{j=1}^{k-1}(1-jt)}.\eeq
\end{prop}

\begin{proof}
Lemma~\ref{lemma:1233123421} translates into the following equations
for the generating functions $R_i$: \bea\nn R_1(t,u,v)&=&tuv
R_2(tv,1,1)+tv
R_3(tv,u,1)+\frac{tv}{v-1}[R_1(t,u,v)-R_1(t,uv,1)], \\
R_2(t,u,v)&=&\frac{tv}{1-t}, \nn\\ R_3(t,u,v)&=&tuv R_1(t,u,1),\label{eqR3} \\
R_4(t,u,v)&=&\frac{tuv}{v-1}[R_1(t,uv,1)-vR_1(t,u,1)].\label{eqR4}\eea
Combining them we get an equation involving only $R_1$:
$$R_1(t,u,v)=\frac{t^2uv^2}{1-tv}+t^2uv^2
R_1(tv,u,1)+\frac{tv}{v-1}[R_1(t,u,v)-R_1(t,uv,1)].$$ If we collect
on one side the terms with $R_1(t,u,v)$, the kernel of the equation
is $1-\frac{tv}{v-1}$. Introducing a new variable $w=uv$ and
canceling the kernel with $v=\frac{1}{1-t}$, we obtain an expression
involving $R_1(t,w,1)$ and $R_1(\frac{t}{1-t},(1-t)w,1)$, which can
be simplified to
$$R_1(t,w,1)=tw^2\left[\frac{1}{1-2t}+\frac{1}{1-t}R_1\left(\frac{t}{1-t},(1-t)w,1\right)\right].$$
By iterated application of this formula, \bea \nn
R_1(t,u,1)&=&t^2u\left(\frac{1}{1-2t}+\frac{t^2u}{1-t}\left(\frac{1}{1-3t}+\frac{t^2u}{1-2t}
\left(\frac{1}{1-4t}+\frac{t^2u}{1-3t}\left(\frac{1}{1-5t}+\cdots\right)\right)\right)\right)\\
\nn
&=&\frac{t^2u}{1-2t}+\frac{(t^2u)^2}{(1-t)(1-3t)}+\frac{(t^2u)^3}{(1-t)(1-2t)(1-4t)}+\frac{(t^2u)^4}{(1-t)(1-2t)(1-3t)(1-5t)}+\cdots\\
\nn
&=&\sum_{k\ge1}\dfrac{t^{2k}u^k}{(1-(k+1)t)\prod_{j=1}^{k-1}(1-jt)}.\eea
Equation (\ref{eqR3}) gives now an expression for $R_3(t,u,1)$, and
(\ref{eqR4}) implies that
$$R_4(t,u,v)=\sum_{k\ge1}\dfrac{t^{2k+1}u^{k+1}(v^2+v^3+\cdots+v^k)}
{(1-(k+1)t)\prod_{j=1}^{k-1}(1-jt)}.$$ Adding up the four generating
functions $R(t,u,1)=R_1(t,u,1)+R_2(t,u,1)+R_3(t,u,1)+R_4(t,u,1)$ we
get (\ref{eq:R}).
\end{proof}

The first coefficients of $R(t,1,1)$, which are the values of
$|\S_n(\{1\mn23,3\mn12,34\mn21\})|$ for $n=1,2,\ldots$, are
$1,2,4,8,19,47,125,\ldots$. This sequence does not appear
in~\cite{Slo} at the moment.

\subsection{$\{1\mn23,34\mn21\}$-avoiding permutations}
\label{sec:threelabsub2}

The derivation of the generating function for this class of
permutations is very similar to the previous subsection. The labels
that we associate to a permutation are again $(s,r,n)$. The proof of
the next lemma is analogous to that of Lemma~\ref{lemma:1233123421}.

\begin{lemma}\label{lemma:1233421} The rightward generating tree for $\{1\mn23,34\mn21\}$-avoiding
permutations is specified by the following succession rule on the
labels: \bce\bt{l}
$(0,1,1)$ \\
$(s,r,n)\longrightarrow\begin{cases} (s+1,1,n+1)\ (s+1,2,n+1)\
\cdots\ (s+1,s,n+1) & \\ \hspace{2cm} (s,s+1,n+1)\ (s,s+2,n+1)\
\cdots\ (s,r,n+1) \quad & \mbox{if } s<r\neq1,
\\ (0,1,n+1)\ (1,2,n+1)\ (1,3,n+1)\ \cdots \ (1,n+1,n+1) & \mbox{if } (s,r)=(0,1),
\\ (s+1,2,n+1)\ (s+1,3,n+1)\ \cdots\ (s+1,s,n+1) & \\
\hspace{2cm} (s,s+1,n+1) \ (s,s+2,n+1)\ \cdots\
 (s,n+1,n+1) & \mbox{if } s>r=1,
\\ \emptyset & \mbox{if } s>r>1.
\end{cases}$
\et\ece
\end{lemma}

Divide the set $\Theta$ of values that the pair $(s,r)$ can take
into four disjoint sets $\Theta_i$, $i=1,2,3,4$ as before, and let
$$T_i(t,u,v):=\sum_{n\ge1}\ \underset{\mathrm{with\
}(s(\pi),r(\pi))\in\Theta_i}{\sum_{\pi\in\S_n(1\mn23,34\mn21)}}
u^{s(\pi)} v^{r(\pi)} \ t^n$$ and
$T(t,u,v)=T_1(t,u,v)+T_2(t,u,v)+T_3(t,u,v)+T_4(t,u,v)$.

\begin{prop}
The generating function for $\{1\mn23,34\mn21\}$-avoiding
permutations where $u$ marks the parameter $s$ defined above is
\beq\label{eq:T}T(t,u,1)=\sum_{k\ge0}\dfrac{t^{k+1}u^k(1+ktu)}{(1+tu)^k(1-kt)(1-(k+1)t)}.\eeq
\end{prop}

\begin{proof}
The equations that follow from Lemma~\ref{lemma:1233421} are now
\bea\nn
T_1(t,u,v)&=&\frac{tuv}{v-1}[T_2(tv,1,1)-T_2(t,1,1)]+\frac{tv}{v-1}
[vT_3(tv,u,1)-T_3(t,uv,1)]\\ \nn &&+\frac{tv}{v-1}[T_1(t,u,v)-T_1(t,uv,1)], \\
T_2(t,u,v)&=&\frac{tv}{1-t}, \nn\\ T_3(t,u,v)&=&tuv T_1(t,u,1),\label{eqT3} \\
T_4(t,u,v)&=&\frac{tuv}{v-1}[T_1(t,uv,1)+T_3(t,uv,1)-vT_1(t,u,1)-vT_3(t,u,1)].\label{eqT4}\eea
From them we can get an equation involving only $T_1$:
$$T_1(t,u,v)=\frac{t^2uv^2}{v-1}\left(\frac{v}{1-tv}-\frac{1}{1-t}\right)+\frac{t^2uv^2}{v-1}
[vT_1(tv,u,1)-T_1(t,uv,1)]+\frac{tv}{v-1}[T_1(t,u,v)-T_1(t,uv,1)].$$
Letting $w=uv$ and canceling the kernel with $v=\frac{1}{1-t}$, we
get that
$$T_1(t,w,1)=\frac{tw}{(1+tw)(1-t)}\left[\frac{t}{1-2t}+T_1\left(\frac{t}{1-t},(1-t)w,1\right)\right].$$
Iterating this formula, we see that
$$T_1(t,u,1)=\sum_{k\ge1}\dfrac{t^{k+1}u^k}{(1+tu)^k(1-kt)(1-(k+1)t)}.$$
Using (\ref{eqT3}) and (\ref{eqT4}) we get expressions for
$T_3(t,u,1)$ and $T_4(t,u,v)$. Finally, the sum
$T(t,u,1)=T_1(t,u,1)+T_2(t,u,1)+T_3(t,u,1)+T_4(t,u,1)$ gives the
formula~(\ref{eq:T}).
\end{proof}

The first coefficients of $T(t,1,1)$ are
$1,2,5,14,42,138,492,\ldots$, which teaches us not to judge a
sequence by looking only at its first five terms. This sequence
gives the number of $\{1\mn23,34\mn21\}$-avoiding permutations of
size $n=1,2,\ldots$, and does not currently appear in~\cite{Slo}.

\section{Concluding remarks}\label{sec:concl}

The main results in the paper have been obtained by constructing
rightward generating trees with up to three labels for several
families of pattern-avoiding permutations, and solving the
functional equations for the generating functions that the
succession rule produces. This is a useful method for enumerating
permutations avoiding generalized patterns. There is nothing special
about the sets of patterns studied in this paper, except that this
method happens to work out nicely on them.

We expect that this technique of rightward generating trees with
several labels, together with the kernel method and other ad-hoc
tools for solving the functional equations that are obtained, will
lead to many more enumerative results for classes of permutations
avoiding generalized patterns.

\subsection*{Acknowledgements}

I am grateful to Mireille Bousquet-M\'elou for many helpful ideas
that have made this paper possible, and to two anonymous referees
for useful suggestions to improve its presentation.


\begin{thebibliography}{99}

\bibitem{BBDFGG} C. Banderier, M. Bousquet-M\'elou, A. Denise, P.
Flajolet, D. Gardy, D. Gouyou-Beauchamps, Generating functions of
generating trees, {\it Discrete Math.} 246 (2002), 29-–55.

\bibitem{BF} C. Banderier and P. Flajolet, Basic analytic
combinatorics of directed lattice paths, {\it Theoret. Comput. Sci.}
281 (2002), 37–-80.

\bibitem{BFP} A. Bernini, L. Ferrari, R. Pinzani, Enumerating
permutations avoiding three Babson-Steingr\'{\i}msson patterns, {\it
Ann. Combin.} 9 (2005), 137--162.

\bibitem{B-M} M. Bousquet-M\'elou, Four classes of pattern-avoiding permutations under one roof: generating trees
with two labels, {\it Electron. J. Combin.} 9 (2003), \#R19.

\bibitem{BP} M. Bousquet-M\'elou and Marko Petkov\v{s}ek, Linear recurrences
with constant coefficients: the multivariate case, {\it Discrete
Math.} 225 (2000), 51-–75.


\bibitem{BEM} A. Burstein, S. Elizalde, T. Mansour, Restricted Dumont
permutations, Dyck paths, and noncrossing partitions, {\it Discrete
Math.} 306 (2006), 2851--2869.

\bibitem{Ca04} D. Callan, Two bijections for Dyck path parameters,
preprint, arxiv:math.CO/0406381v2.

\bibitem{C}
A. Claesson, Generalised pattern avoidance, {\it Europ. J. Combin.}
22 (2001), 961--973.

\bibitem{ClaMan}
A. Claesson, T. Mansour, Enumerating permutations avoiding a pair of
Babson-Steingr\'{\i}msson patterns, {\it Ars Combinatorica} 77
(2005).

\bibitem{Eli06} S. Elizalde, Asymptotic enumeration of permutations avoiding generalized
patterns, {\it  Adv. in Appl. Math.}  36  (2006),  138--155.

\bibitem{EliMan} S. Elizalde, T. Mansour, Restricted Motzkin permutations, Motzkin paths,
continued fractions, and Chebyshev polynomials, {\it Discrete Math.}
305 (2005), 170--189.

\bibitem{EliNoy} S. Elizalde, M. Noy, Consecutive subwords in permutations,
{\it Adv. in Appl. Math.} 30 (2003), 110--125.

\bibitem{FirMan} G. Firro, T. Mansour, Three-letter-pattern-avoiding permutations and functional
equations, {\it Electron. J. Combin.} 13 (2006), \#R51.

\bibitem{Kra} C. Krattenthaler, Permutations with restricted
patterns and Dyck paths, {\it Adv. in Appl. Math.} 27 (2001),
510--530.

\bibitem{Man} T. Mansour, Restricted $1\mn3\mn2$ permutations and generalized
patterns, {\it Ann. Combin.} 6 (2002), 65-76.

\bibitem{Slo} N.J.A. Sloane, S. Plouffe, The Encyclopedia of Integer Sequences,
Academic Press, San Diego, 1995,
\texttt{http://www.research.att.com/$\sim$njas/sequences}.

\bibitem{West} J. West, Generating trees and the Catalan and Schr\"oder numbers,
{\it Discrete Math.} 146 (1995), 247--262.

\bibitem{Wes96} J. West, Generating trees and forbidden
subsequences, {\it Discrete Math.} 157 (1996), 363--374.

\end{thebibliography}
\end{document}